\documentclass[oneside,a4paper,reqno]{amsart}
\usepackage{pdfsync, enumerate}
\usepackage{stmaryrd}
\usepackage{mathrsfs}
\usepackage{multicol}
\usepackage{amsmath, amsthm, amscd, amssymb, latexsym, eucal}
\usepackage[all]{xy}
\def\serieslogo@{} \def\@setcopyright{} \makeatother

\usepackage{multienum}

\usepackage[colorinlistoftodos]{todonotes}

\usepackage{hyperref}
\usepackage{color}
\usepackage{cite}

\makeatletter
\renewcommand*\env@matrix[1][c]{\hskip -\arraycolsep
  \let\@ifnextchar\new@ifnextchar
  \array{*\c@MaxMatrixCols #1}}
\makeatother

\usepackage{color}

\numberwithin{equation}{section}
\newtheorem{thm}{Theorem}[section]
\newtheorem*{main-thm}{Main Theorem}
\newtheorem*{Auslander-thm}{Auslander's Theorem}

\newtheorem{cor}[thm]{Corollary}
\newtheorem{lem}[thm]{Lemma}
\newtheorem{prop}[thm]{Proposition}

\theoremstyle{definition}
\newtheorem{defn}[thm]{Definition}

\newtheorem{exam}[thm]{Example}

\newtheorem*{problem}{Problem}

\newcommand{\lxr}{\longrightarrow}

\newcommand{\A}{\mathscr A}
\newcommand{\B}{\mathscr B}
\newcommand{\C}{\mathscr C}

\newcommand{\mr}{\mathsf{r}}
\newcommand{\mq}{\mathsf{q}}
\newcommand{\mi}{\mathsf{i}}

\newcommand{\ml}{\mathsf{l}}
\newcommand{\me}{\mathsf{e}}

\newcommand{\map}{\mathsf{p}}

\newcommand{\mD}{\mathsf{D}}

 \DeclareMathOperator*{\opp}{\mathsf{op}}

\DeclareMathOperator*{\Ker}{\mathsf{Ker}}
 \DeclareMathOperator*{\Image}{\mathsf{Im}}

 \DeclareMathOperator{\pd}{\mathsf{pd}}
\DeclareMathOperator*{\id}{\mathsf{id}}

 \DeclareMathOperator*{\smod}{\mathsf{mod}-\!}

\DeclareMathOperator*{\proj}{\mathsf{proj}}

\DeclareMathOperator*{\Proj}{\mathsf{Proj}}

 \DeclareMathOperator{\Hom}{\mathsf{Hom}}
 \DeclareMathOperator{\HH}{\mathsf{HH}}

\DeclareMathOperator{\Ext}{\mathsf{Ext}}

\newcommand{\evl}{{\operatorname{env}}}

\newcommand{\perf}{{\operatorname{perf}}}
\newcommand{\rad}{\mathsf{rad}}

\newcommand{\iden}{\operatorname{Id}\nolimits}

\newcommand{\fg}[1]{\textup{\textsf{Fg#1}}}
\newcommand{\Dsg}{\mathsf{D}_\mathsf{sg}}

\newsavebox{\proofbox}
\savebox{\proofbox}{\begin{picture}(7,7)%
  \put(0,0){\framebox(7,7){}}\end{picture}}

\usepackage{latexsym}
\usepackage{pstricks}
\usepackage{comment}

\begin{document}

\title[Arrow removal]{Homological invariants of the arrow removal operation}

\author[Erdmann]{Karin Erdmann}
\address{Karin Erdmann\\
Mathematical Institute\\ 
24--29 St.\ Giles\\
Oxford OX1 3LB\\ 
England}
\email{erdmann@maths.ox.ac.uk}

\author[Psaroudakis]{Chrysostomos Psaroudakis}
\address{Department of Mathematics, Aristotle University of Thessaloniki, Thessaloniki 54124, Greece}
\email{chpsaroud@math.auth.gr}

\author[Solberg]{\O yvind Solberg}
\address{Department of Mathematical Sciences\\
NTNU\\
N-7491 Trondheim, Norway }
\email{oyvind.solberg@math.ntnu.no}

\date{\today}
 


\begin{abstract}
  In this paper we show that Gorensteinness, singularity categories
  and the finite generation condition \fg{} for the Hochschild
  cohomology are invariants under the arrow removal operation for a
  finite dimensional algebra.
\end{abstract}

\maketitle

\setcounter{tocdepth}{1} \tableofcontents

\section{Introduction and the main result}

In \cite{GPS} the arrow removal operation on quotients of path
algebras was investigated with respect to the finitistic dimension
conjecture. The idea was to remove those arrows that don't contribute
to the finitistic dimension. This new technique gave us a successful
reduction method for actual computing the finistic dimension in many
examples. Our aim in this paper is to investigate further this
operation, and in particular study the class of {\em arrow removal
  algebras} with respect to various homological invariants. The
guiding problem can be formulated as follows$\colon$

\begin{problem}
  How does the arrow removal behave with respect to Gorensteinness,
  singularity categories and the finite generation condition \fg{} for
  the Hochschild cohomology.
\end{problem}
  See Section \ref{section:4} for the definition of
  Gorensteinness \cite{AuslanderReiten2} and the singularity category $\Dsg(\Lambda)$ of a
  finite dimensional algebra $\Lambda$ \cite{Buchweitz}, and Section \ref{section:5}
  for the definition of \fg{} \cite{EHSST, Solberg:contemp}.

In our main result we prove that the above three homological
invariants remain the same under the arrow removal operation for an
admissible path algebra over a field. Before we state our main result
we briefly recall the arrow removal operation.

Let $\Lambda$ be an admissible quotient $kQ/I$ of a path algebra $kQ$
over a field $k$.  Consider an arrow $a$ in $Q$ such that $a$ does not
occur in a minimal generating set of $I$. Then the quotient algebra
$\Gamma=\Lambda/\langle \overline{a}\rangle$ is called an {\em
  arrow removal algebra} of $\Lambda$.  This new algebra can be
explicitly described as a trivial extension. More precisely, it has
been proved in \cite[Theorem~A]{GPS} that the arrow
$a\colon v_e\to v_f$ in $Q$ does not occur in a set of minimal
generators of $I$ in $kQ$ if and only if
$\Lambda\cong \Gamma\ltimes P$, where $P= \Gamma e\otimes_kf\Gamma$
with $\Hom_\Gamma(e\Gamma, f\Gamma)=(0)$. Here $e$ is a trivial path
in $kQ$ and $v_e$ denotes the corrsponding vertex in $Q$.  We can
  also consider arrow removal for any finite number of arrows, but for
  simplicity we only review our results in the one arrow case in the
  introduction.

The module category of $\Lambda$ can be described using the context of
cleft extensions of abelian categories in the sense of Beligiannis
\cite{Bel:Cleft}. This means that there are functors
$\mi\colon \smod\Gamma\to \smod\Lambda$ and
$\me\colon \smod\Lambda\to \smod\Gamma$, induced by the natural
surjection map $\Gamma\to \Lambda$ and by the inclusion map
$\Gamma\to \Lambda$ respectively, the functor $\me$ is faithful exact
and the composition $\me\circ \mi$ is equivalent to the identity
functor on $\smod{\Gamma}$. The arrow removal operation gives even
more homological properties on this cleft extension.

The notion of eventually homological isomorphism was introduced in
\cite{PSS}.  Recall that a functor $F\colon \B\to \C$ between abelian
categories is called an {\em eventually homological isomorphism}, if
there is an integer $t$ such that for every $j>t$ there is an
isomorphism $\Ext_{\B}^j(X,Y)\cong \Ext_{\C}^j(FX,FY)$ for all objects
$X, Y'\in \B$. Given the smallest such $t$, we call the functor a
  \emph{$t$-eventually homological isomorphism}.  The latter notion
was used in comparing the algebras $\Lambda$ and $e\Lambda e$
for an idempotent $e$ in $\Lambda$, with respect to Gorensteinness,
singularity categories and the finite generation condition \fg{} for
the Hochschild cohomology (\!\!\cite[Main Theorem]{PSS}).  We mention
that this comparison theorem was achieved via recollements of abelian
categories.

We summarize below our main results in the simplified setting of a one
arrow removal.

\begin{main-thm}
  Let $\Lambda = kQ/I$ be an admissible quotient of a path algebra
  $kQ$ over a field $k$ and let
  $\Gamma=\Lambda/\langle \overline{a}\rangle$ an arrow removal of
  $\Lambda$ for an arrow $a$ in $Q$. Then the following hold.
\begin{enumerate}[\rm(i)]
\item The functor $\me\colon\smod\Lambda\to \smod\Gamma$ is a
  $1$-eventually homological isomorphism.

\item $\Lambda$ is Gorenstein if and only if $\Gamma$ is Gorenstein.

\item The functor $\me\colon \Dsg(\Lambda) \to \Dsg(\Gamma)$ is a
  singular equivalence. 

\item $\Lambda$ satisfies
  \fg{} if and only if $\Gamma$ satisfies \fg{}.
\end{enumerate}
\end{main-thm}

We remark that the arrow removal operation has been also considered in
\cite{CLMS}. They mainly worked on the converse process, i.e.\ add
arrows to a path algebra, and they described the Hochschild
(co)homology using different techniques.


We end the introduction with a short description of the contents of
the paper section by section.  In Section \ref{sectioncleftextabcat}
we review relevant results on cleft extensions.  As part of a cleft
extension between abelian categories $\mathcal{A}, \mathcal{B}$ there
is a functor $\me\colon \mathcal{A}\to \mathcal{B}$ which is faithful
and exact. One would like this to be an eventually homological
isomorphism.  In Section \ref{section:3}, we show this is the case
under certain conditions.  Section \ref{section:4} shows that
Gorensteinness and singularity categories are invariant under arrow
removal. In Section \ref{section:5} we investigate the Fg condition,
and prove that it is invariant under arrow removal.



\section{Cleft extensions and arrow removals}
\label{sectioncleftextabcat}

We start this section by recalling and reviewing some results about
cleft extensions of abelian categories from \cite{Bel:Cleft,GPS} that we need in the sequel.

\subsection{Cleft extensions}
We first recall the definition of cleft extensions of abelian
categories.

\begin{defn}\textnormal{(\!\!\cite[Definition 2.1]{Bel:Cleft})}
\label{defncleftext}
A {\bf cleft extension} of an abelian category $\B$ is an
abelian category $\A$ together with functors:
\[
\xymatrix@C=0.5cm{
\B \ar[rrr]^{\mi} &&& \A \ar[rrr]^{\me} &&& \B
\ar @/_1.5pc/[lll]_{\ml} } 
\]
henceforth denoted by $(\B,\A, \me, \ml, \mi)$, such that the
following conditions hold:
\begin{enumerate}[\rm(a)]
\item The functor $\me$ is faithful exact.

\item The pair $(\ml,\me)$ is an adjoint pair of functors,
  where we denote the adjunction by
\[
\theta_{B,A}\colon \Hom_\A(\ml(B), A) \simeq \Hom_\B(B,\me(A)).
\]

\item There is a natural isomorphism
  $\varphi\colon \me\mi\lxr \mathsf{Id}_{\B}$ of functors.
\end{enumerate} 
\end{defn}

Denote the unit $\theta_{B,\ml(B)}(1_{\ml(B)})$ and the counit
$\theta^{-1}_{\me(A),A}(1_{\me(A)})$ of the adjoint pair $(\ml,\me)$
by $\nu\colon 1_\B\lxr \me\ml$ and $\mu\colon \ml\me\lxr 1_\A$,
respectively.  The unit and the counit satisfy the relations 
\begin{equation}\label{eq:unit}
1_{\ml(B)} = \mu_{\ml(B)}\ml(\nu_B)
\end{equation}
and
\begin{equation}\label{eq:counit}
1_{\me(A)} = \me(\mu_A)\nu_{\me(A)}
\end{equation}
for all $B$ in $\B$ and $A$ in $\A$.  From \eqref{eq:counit} the
morphism $\me(\mu_A)$ is an (split) epimorphism. Since $\me$ is
faithful exact, it follows that $\mu_A$ is an epimorphism for all $A$
in $\A$. This implies that for all $A$ in $\A$ the following sequence
is exact
\begin{equation}
\label{firstfundamentalsequence}
\xymatrix@C=0.5cm{
0 \ar[rr] && \Ker\mu_A \ar[rr]^{} && \ml\me(A) \ar[rr]^{\mu_A} && A
\ar[rr] && 0  }
\end{equation}

The next result collects some basic properties of a cleft extension
which basically follows from Definition~\ref{defncleftext} and are
discussed in \cite{Bel:Cleft}. For a detailed proof
the reader is referred to \cite[Lemma~2.2]{GPS}.

\begin{lem}
\label{basicproperties}
Let $\A$ be a cleft extension of $\B$. Then the following hold.
\begin{enumerate}[\rm(i)]
\item The functor $\me\colon \A\lxr \B$ is essentially surjective. 

\item The functor $\mi\colon \B \lxr \A$ is fully faithful and exact.

\item The functor $\ml\colon \B\lxr \A$ is faithful and preserves
  projective objects.

\item There is a functor $\mq\colon \A\lxr \B$ such that $(\mq,\mi)$
  is an adjoint pair.

\item There is a natural isomorphism $\mq\ml\simeq \mathsf{Id}_{\B}$ of
  functors. 
\end{enumerate}
\end{lem}

A cleft extension $(\B,\A, \me, \ml, \mi)$ is equipped with three
  additional functors that are crucial in our investigations.  We saw
  in \eqref{firstfundamentalsequence} that there is a short exact
  sequence
\[\xymatrix@C=0.5cm{
0 \ar[rr] && \Ker\mu_A \ar[rr]^{} && \ml\me(A) \ar[rr]^{\mu_A} && A
\ar[rr] && 0  }\]
for all $A$ in $\A$.  The assignment $A\mapsto \Ker\mu_A$ defines an
endofunctor $G\colon \A\lxr \A$ and therefore also an exact sequence
of endofunctors on $\A$ 
\begin{equation}\label{eq:G-functor}
\xymatrix@C=0.5cm{
0 \ar[rr] && G \ar[rr]^{} && \ml\me \ar[rr]^{\mu_{(-)}} && \mathsf{Id}_\A
\ar[rr] && 0  }
\end{equation}
Precompose the above exact sequence of functors with the functor
$\mi\colon \B \to \A$, and we obtain an  exact sequence of functors
\[\xymatrix@C=0.5cm{
0 \ar[rr] && G\mi \ar[rr]^{} && \ml\me\mi \ar[rr]^{\mu_{\mi(-)}} &&
\mi \ar[rr] && 0  }\]
Denote the functor $G\mi \colon \B\to \A$ by $H$ and view
$\varphi\colon \me\mi\to \mathsf{Id}_{\B}$ as an identification.  Then
we have the exact sequence of funtors 
\begin{equation}\label{eq:H-functor}
\xymatrix@C=0.5cm{
0 \ar[rr] && H \ar[rr]^{} && \ml \ar[rr]^{\mu_{\mi(-)}} &&
\mi \ar[rr] && 0  }
\end{equation}
Postcompose the above exact sequence with the functor $\me\colon
\A\to\B$ and obtain the exact sequence 
\[0\to \me H \to \me\ml \xrightarrow{\me(\mu_{\mi(-)})} \me\mi \to 0.\]
Again, viewing $\varphi\colon \me\mi\to \mathsf{Id}_\B$ as an
identification and denote the endofunctor $\me H$ on $\B$ by $F$.
Then we obtain an exact sequence of endofunctors on $\B$ 
\begin{equation}
\label{splitsequencegivesF}
\xymatrix@C=0.5cm{
0 \ar[rr] && F \ar[rr]^{} && \me\ml \ar[rr]^{ \ \ \me(\mu_{\mi(-)})}
&& \mathsf{Id}_\B \ar[rr] && 0  }
\end{equation}
The following lemma is an immediate consequence of \eqref{eq:counit}. 
\begin{lem}
\label{lemsplitexactseq}
Let $(\B,\A,\me, \ml, \mi)$ be a cleft extension of abelian
categories. Then the exact sequence $(\ref{splitsequencegivesF})$
splits.
\end{lem}

Another fact on cleft extensions we use later, is the following result (see
  \cite[Lemma 2.4]{GPS}).
  
\begin{lem}
\label{lemnilpotentfunctors}
Let $(\B,\A,\me, \ml, \mi)$ be a cleft extension of abelian categories.
The following statements hold.
\begin{enumerate}[\rm(i)]
\item For any $n\geq 1$, there is a natural isomorphism $\me G^n\simeq F^n\me$.

\item Let $n\geq 1$. Then $F^n=0$ if and only if $G^n=0$.
\end{enumerate}
\end{lem}

In Sections \ref{section:3} and \ref{section:5} the following assumption on a cleft
extension $(\B,\A,\me, \ml, \mi)$ of abelian categories shall be of
importance.
\begin{equation}
\label{assumptionscleftext}
\text{The functor} \ \ml \ \text{is exact and the functor} \ \me \ \text{preserves projectives.}
\end{equation}

\subsection{Cleft extensions arising from arrow removals}\label{subsec:cleftextfromarrowremoval}
Let $\Lambda=kQ/I$ be an admissible quotient of a path algebra $kQ$
over a field $k$.  Suppose that there is a set of arrows
$a_i\colon v_{e_i}\to v_{f_i}$ in $Q$ for $i=1,2,\ldots,t$ which do
not occur in a set of minimal generators of $I$ in $kQ$ and
$\Hom_\Lambda( e_i\Lambda, f_j\Lambda ) = 0$ for all $i$ and $j$ in
$\{1,2,\ldots,t\}$.  Let
$\Gamma = \Lambda/\Lambda\{\overline{a}_i\}_{i=1}^t\Lambda$.  Recall
that the natural projection $\pi\colon \Lambda \to \Gamma$ or just the
pair $\Lambda$ and
$\Gamma = \Lambda/\Lambda\{\overline{a}_i\}_{i=1}^t\Lambda$ is an
\emph{arrow removal}. The following result from \cite{GPS} shows that
the arrow removal operation induces a cleft extension between the
corresponding module categories, i.e.\ $\smod\Lambda$ is a cleft
extension of $\smod\Gamma$, with certain homological properties.

\begin{thm}\label{cleftextarrowrem} \textnormal{(\!\! For (i)
    \cite[Corollary~4.3, Proposition~4.6]{GPS}, and for (ii)
    Proposition~4.6)}
\label{propremoveiscleft}
Let $\Lambda = kQ/I$ be a quotient path algebra as above.  
\begin{enumerate}[\rm(i)]
\item A set of arrows $a_i\colon v_{e_i}\to v_{f_i}$ in $Q$ for
  $i=1,2,\ldots,t$ which do not occur in a set of minimal generators
  of $I$ in $kQ$ and $\Hom_\Lambda( e_i\Lambda, f_j\Lambda ) = 0$ for
  all $i$ and $j$ in $\{1,2,\ldots,t\}$ if and only if $\Lambda$ is
  isomorphic to the trivial extension $\Gamma\ltimes P$, where
  $\Gamma = \Lambda/\Lambda \{\overline{a}_i\}_{i=1}^t \Lambda$ and
  $P = \oplus_{i=1}^t \Gamma e_i\otimes_kf_i\Gamma$ with
  $\Hom_\Gamma(e_i\Gamma,f_j\Gamma)=0$ for all $i,j = 1,2,\ldots,t$.

\item Suppose that there are arrows $a_i\colon v_{e_i}\lxr v_{f_i}$ in
  $Q$ for $i=1,2,\ldots,t$ which do not occur in a set of minimal
  generators of $I$ in $kQ$ and
  $\Hom_\Lambda(e_i\Lambda, f_j\Lambda ) = 0$ for all $i$ and $j$ in
  $\{1,2,\ldots,t\}$.  Let $\Gamma = \Lambda/\Lambda
\{\overline{a}_i\}_{i=1}^t \Lambda$.   Then the tuple
  $(\smod\Gamma, \smod\Lambda, \me, \ml, \mi)\colon$
\begin{equation}
\label{digramwithendofunctors}
\xymatrix@C=0.5cm{
\smod\Gamma  \ar@(ur,ul)_{F} \ar[rrr]^{\mi = \Hom_\Gamma({_\Lambda\Gamma}_\Gamma,-)} &&& \smod\Lambda
\ar@(ur,ul)_{G} \ar[rrr]^{\me = \Hom_\Lambda({_\Gamma\Lambda}_\Lambda,-)} 
\ar @/_1.5pc/[lll]_{\mq = {-\otimes_\Lambda  {_\Lambda\Gamma}_\Gamma}}  \ar
 @/^1.5pc/[lll]^{\map = \Hom_\Lambda({_\Gamma\Gamma}_\Lambda,-)} &&& \smod\Gamma \ar@(ul,ur)^{F}
\ar @/_1.5pc/[lll]_{\ml = {-\otimes_\Gamma {_\Gamma\Lambda}_\Lambda}} \ar
 @/^1.5pc/[lll]^{\mr = \Hom_\Gamma({_\Lambda\Lambda}_\Gamma,-)} }
\end{equation}
is a cleft extension satisfying the following conditions, where
  $F$ and $G$
are as in \eqref{splitsequencegivesF} and \eqref{eq:G-functor}: 
 \begin{enumerate}[\rm(a)]
\item $\me$ is faithful exact,
\item $(\ml,\me)$ is an adjoint pair of functors,
\item $\me\mi \simeq 1_{\smod\Gamma}$,
\item $\ml$ and $\mr$ are exact functors,
\item $\me$ preserves projectives, 
\item $\Image F \subseteq \proj(\Gamma)$ and $\Image G\subseteq
  \Proj(\Lambda)$, 
\item $F^2 = 0$. 
\end{enumerate}
\end{enumerate}
\end{thm}
From the proof of \cite[Proposition 4.6 (iv)]{GPS} we have the
following description of the functors $F$ and $F^{\opp}$ (when we
consider left modules).

\begin{lem}\label{lem:F-func-prop}
  Let $\Lambda=kQ/I$ be an admissible quotient of a path algebra $kQ$
  over a field $k$.  For a set of arrows
  $a_i\colon v_{e_i}\to v_{f_i}$ in $Q$ for $i=1,2,\ldots,t$ suppose
  that $\Lambda \to \Gamma =
  \Lambda/\Lambda\{\overline{a}\}_{i=1}^t\Lambda$ is an arrow
  removal. Then 
\begin{enumerate}[\rm(a)]
\item The endofunctor $F\colon \smod\Gamma\to \smod\Gamma$ is given as 
\[F = -\otimes_\Gamma \Gamma\{\overline{a}_i\}_{i=1}^t\Gamma\colon
  \smod\Gamma\to \smod\Gamma.\]
\item The endofunctor $F^{\opp}\colon \smod\Gamma^{\opp}\to \smod\Gamma^{\opp}$ is given as 
\[F^{\opp} = \Gamma\{\overline{a}_i\}_{i=1}^t\Gamma\otimes_\Gamma - \colon
  \smod\Gamma^{\opp}\to \smod\Gamma^{\opp}.\]
\end{enumerate}
\end{lem}

We remark that the above homological properties were used to show that
the finiteness of the finitistic dimension of $\Lambda$ can be reduced
to the finiteness of the finitistic dimension of the arrow removal
algebra $\Gamma$, see \cite[Theorem~A]{GPS}.  These homological
properties are also used intensively in the sequel of the paper to
show the invariance of Gorensteinness, singularity categories and the
finite generation condition \fg{} for the Hochschild cohomology under
the arrow removal operation. It is interesting that this operation
gives rise to such a powerful cleft extension.

\section{Cleft extensions and eventually homological
  isomorphisms}\label{section:3} 

Let $(\B,\A, \me, \ml, \mi)$ be a cleft extension of abelian
categories.  Then the functor $\me\colon \A\to \B$ is an exact
functor, so it always induces homomorphism of the following Yoneda
rings  
\[
\me\colon \Ext^*_\A(A,A) \to \Ext^*_\B( \me(A),\me(A) )
\]
for all $A$ in $\A$.  Recall from \cite[Section 3]{PSS} that $\me$ is
called an \emph{eventually homological isomorphism} if 
\[\Ext^i_\A(A,A) \simeq \Ext^i_\B( \me(A),\me(A) )\]
for every $i > t$ where $t$ is some positive integer. For the minimal
such $t$, the above isomorphism is called a \emph{$t$-eventually
  homological isomorphism}.  Note that in the definition we do not
require that the isomorphism is induced by the functor $\me$.

In this section we describe one situation where the functor $\me$ of a
cleft extension is an eventually homological isomorphism. We start
with the following result.

\begin{lem}\label{lem:commdiagcleft}
  Let $(\B,\A, \me, \ml, \mi)$ be a cleft extension of abelian
  categories such that condition \textup{\eqref{assumptionscleftext}} is
  satisfied.  Then, for all $i\geq 1$ and for all $A, C\in \A$, the
  following diagram commutes
\[
\xymatrix{
\Ext^i_\A(C,A) \ar[r]^-\me\ar@{=}[d] & \Ext^i_\B(\me(C),\me(A)) \\
\Ext^i_\A(C,A) \ar[r]^-{\mu_C^*} & \Ext^i_\A(\ml\me(C),A)
\ar[u]_{\theta_{\me(C),A}}^{\simeq} }
\]
where the vertical maps are isomorphisms. 
\end{lem}
\begin{proof}
Let $m\geq 1$, and let $f\in\Ext^m_\A(C,A)$ be represented by a
morphism $f\colon \Omega_\A^m(C)\to A$.  Then consider the following exact
commutative diagram 
\[\xymatrix{
0\ar[r] & \Omega_\A^m(C) \ar[r]\ar[d]_f & P_{m-1} \ar[r]\ar[d] & P_{m-2}\ar[r]\ar@{=}[d] &
\cdots \ar[r] & P_0\ar[r]\ar@{=}[d] & C \ar[r]\ar@{=}[d] & 0\\
0\ar[r] & A \ar[r] & M \ar[r] & P_{m-2}\ar[r] &
\cdots\ar[r] & P_0\ar[r] & C \ar[r] & 0
}\]
Apply the exact functor $\me$ to this diagram and obtain the following
diagram
\[\xymatrix@C15pt{
0\ar[r] & \me(\Omega_\A^m(C)) \ar[r]\ar[d]_{\me(f)} & \me(P_{m-1}) \ar[r]\ar[d] & \me(P_{m-2})\ar[r]\ar@{=}[d] &
\cdots \ar[r] & \me(P_0)\ar[r]\ar@{=}[d] & \me(C) \ar[r]\ar@{=}[d] & 0\\
0\ar[r] & \me(A) \ar[r] & \me(M) \ar[r] & \me(P_{m-2})\ar[r] &
\cdots\ar[r] & \me(P_0)\ar[r] & \me(C) \ar[r] & 0
}\]
The lower row represents the image of $f$ under the functor $\me$. 

By the adjunction
$\theta\colon \Hom_\B(\me(\Omega_\A^m(C)),\me(A)) \simeq
\Hom_\A(\ml\me(\Omega_\A^m(C)), A)$ the morphism $\me(f)$ corresponds
to $\theta(\me(f))\colon \ml\me(\Omega_\A^m(C))\to A$.  This last
morphism is equal to the composition of the morphisms
\[
\ml\me(\Omega_\A^m(C))\xrightarrow{\ml\me(f)} \ml\me(A)
  \xrightarrow{\mu_A} A.
\]
In addition we have the following two commutative diagrams
\[
\xymatrix{
0\ar[r] & \ml\me(A) \ar[r]\ar[d]^{\mu_A} & \ml\me(M) \ar[r]\ar[d]^s &
\ml\me(\Omega_\A^{m-1}(C))\ar[r]\ar@{=}[d] & 0\\
0\ar[r] & A \ar[r]\ar@{=}[d] & E \ar[r]\ar[d]^t &
\ml\me(\Omega_\A^{m-1}(C))\ar[r]\ar[d]^{\mu_{\Omega_\A^{m-1}(C)}} & 0\\
0\ar[r] & A \ar[r] & M \ar[r] & \Omega_\A^{m-1}(C)\ar[r] & 0
}
\]
and 
\[
\xymatrix{
0\ar[r] & \ml\me(\Omega_\A^1(C)) \ar[r]\ar[d]^{\mu_{\Omega_\A^1(C)}} & \ml\me(P_0) \ar[r]\ar[d]^{s'} &
\ml\me(C)\ar[r]\ar@{=}[d] & 0\\
0\ar[r] & \Omega_\A^1(C) \ar[r]\ar@{=}[d] & E' \ar[r]\ar[d]^{t'} &
\ml\me(C)\ar[r]\ar[d]^{\mu_{C}} & 0\\
0\ar[r] & \Omega_\A^1(C) \ar[r] & P_0 \ar[r] & C\ar[r] & 0
}
\]
with $ts = \mu_{M}$ and $t's' =
\mu_{P_0}$. Using these commutative diagrams we can
construct the following commutative diagram 
\[
\xymatrix@C11pt{
0\ar[r] & \ml\me(\Omega_\A^m(C)) \ar[r]\ar[d]_{\ml\me(f)} &
\ml\me(P_{m-1}) \ar[r]\ar[d] & \ml\me(P_{m-2})\ar[r]\ar@{=}[d] & 
\cdots \ar[r] & \ml\me(P_1)\ar[r]\ar@{=}[d] & \ml\me(P_0)\ar[r]\ar@{=}[d] & \ml\me(C) \ar[r]\ar@{=}[d] & 0\\
0\ar[r] & \ml\me(A) \ar[r]\ar[d]_{\mu_A} & \ml\me(M) \ar[r]\ar[d]^s & \ml\me(P_{m-2})\ar[r]\ar@{=}[d] &
\cdots\ar[r] & \ml\me(P_1)\ar[r]\ar@{=}[d]
&\ml\me(P_0)\ar[r]\ar@{=}[d] & \ml\me(C) \ar[r]\ar@{=}[d] & 0\\
0\ar[r] & A \ar[r]\ar@{=}[d] & E \ar[r]\ar[d]^t & \ml\me(P_{m-2})\ar[r]\ar[d]^{\mu_{P_{m-2}}} &
\cdots\ar[r] & \ml\me(P_1)\ar[r]\ar[d]^{\mu_{P_1}}
&\ml\me(P_0)\ar[r]\ar[d]^{s'} & \ml\me(C) \ar[r]\ar@{=}[d] & 0\\
0\ar[r] & A \ar[r]\ar@{=}[d] & M \ar[r]\ar@{=}[d] &
P_{m-2}\ar[r]\ar@{=}[d] &\cdots\ar[r]
& P_1\ar[r]\ar@{=}[d] & E'\ar[r]\ar[d]^{t'} & \ml\me(C)\ar[r]\ar[d]^{\mu_C} & 0\\
0\ar[r] & A\ar[r] & M\ar[r] & P_{m-2}\ar[r] & \cdots\ar[r] & P_1\ar[r]
& P_0\ar[r] & C\ar[r] & 0
}
\]
The third row in the above diagram corresponds to the
$\theta^{-1}_{\me(C),A}(\me(f))$, and the fourth row in the above
diagram corresponds to the image of $f$ under the map $\mu_C^*$.  It
follows that the diagram in the statement is commutative. 

Finally, using the adjunction $(\ml, \me)$ and since both functors are
exact and preserve projectives, it follows immediately that the map
$\theta_{\me(C),A}$ is an isomorphism.
\end{proof}
If the map induced by $\mu^*_C$ in Lemma \ref{lem:commdiagcleft} is an isomorphism
for all $C$ and all $i\gg 0$, it would follow that the functor $\me$
is an eventually homological isomorphism.  Using a homological
condition on the functor $G$ (see \eqref{eq:G-functor}) the next
result describes a situation when $\mu^*_C$ induces such an isomorphism.

\begin{thm}\label{thm:cleftabelianeventuallyisom}
  Let $(\B,\A, \me, \ml, \mi)$ be a cleft extension of abelian
  categories satisfying condition \textup{\eqref{assumptionscleftext}}.  Assume
  that 
\[\sup\{\pd_\A G(A) \mid A\in\A\} \leq n_\A\]
for some integer $n_\A$.  Then the functor $\me\colon \A\to \B$ is an
$n_\A + 1$-eventually homological isomorphism.
\end{thm}
\begin{proof}  
Using the commutative diagram in Lemma \ref{lem:commdiagcleft} and the
long exact sequence induced from the exact sequence
\[0\to G(C)\to \ml\me(C) \to C \to 0\]
applying the functor $\Ext^*_\A(-,A)$, in fact the functor $\me$
induces an isomorphism between $\Ext^i_\A(C,A)$ and
$\Ext^i_\B(\me(C),\me(A))$ for $i > n_\A + 1$.  The claim follows from
this.
\end{proof}

Applying Theorem~\ref{thm:cleftabelianeventuallyisom} to the cleft
extension of the arrow removal, see Theorem~\ref{cleftextarrowrem}, we
get the following consequence. This result constitutes part (i) of the
Main Theorem presented in the Introduction.

\begin{cor}
\label{corarrowremovalehi}
Let $\Lambda = kQ/I$ be an admissible quotient of a path algebra $kQ$
over a field $k$, and assume that
$\Gamma=\Lambda/\langle \{\overline{a}_i\}_{i=1}^t\rangle$ is an arrow
removal of $\Lambda$ for the arrows $\{a_i\}_{i=1}^t$ in $Q$.  Then
the functor $\me\colon\smod\Lambda\to \smod\Gamma$ is a
$1$-eventually homological isomorphism.
\end{cor}

\section{Gorenstein algebras and singular
  equivalences}\label{section:4}

In this section we show that Gorensteinness and singularity categories
are invariant under the arrow removal operation.  Recall from
  \cite{AuslanderReiten2} that a finite dimensional algebra $\Lambda$ is
  called \emph{Gorenstein} if $\Lambda$ satisfies
  $\id_{\Lambda^{\opp}} {_\Lambda\Lambda} < \infty$ and
  $\id_{\Lambda} \Lambda_\Lambda < \infty$.  Furthermore,
  recall 
  from \cite{Buchweitz} that the \emph{singularity category}
  $\Dsg(\A)$ of an abelian category $\A$ with enough projectives is
  given by the Verdier quotient $\mD^b(\A)/\mD^\perf(\A)$. Here
  $\mD^\perf(\A)$ denotes the full triangulated subcategory of
  $\mD^b(\A)$ consisting of the perfect objects, i.e.\ complexes
  quasi-isomorphic to bounded complexes with components in $\Proj\A$.

\subsection{Gorenstein algebras}
Let $\Lambda=kQ/I$ be an admissible quotient of a path algebra $kQ$,
and suppose $\{a_i\}_{i=1}^t$ is a set set of arrows in $Q$ such that
$\Gamma = \Lambda/\langle \{\overline{a}_i\}_{i=1}^t\rangle$ is an arrow
removal. The key fact for the invariance of Gorensteiness is that the
functor $\me\colon\smod\Lambda\to \smod\Gamma$ is an eventually
homological isomorphism, as shown in
Corollary~\ref{corarrowremovalehi}. The reason is the following result
which we formulate, for simplicity, for module categories over finite
dimensional algebras.

\begin{thm}\textnormal{(\!\!\cite[Theorem 4.3 (v)]{PSS})}
\label{thmehi14}
Let $T\colon \smod\Lambda\to \smod \Gamma$ be a functor which is
essentially surjective and an eventually homological isomorphism. Then
$\Lambda$ is Gorenstein if and only if $\Gamma$ is Gorenstein.
\end{thm}

We can now show that Gorensteinness is indeed invariant under the
arrow removal operation. In particular, the following result is an
immediate consequence of Corollary~\ref{corarrowremovalehi} and
Theorem~\ref{thmehi14}. This result constitutes part (ii) of the Main
Theorem presented in the Introduction.

\begin{cor}
\label{thmGorenstein}
Let $\Lambda=kQ/I$ be an admissible quotient of a path algebra $kQ$,
and suppose that
$\Gamma=\Lambda/\langle \{\overline{a}_i\}_{i=1}^t\rangle$ is an arrow
removal of $\Lambda$ for the arrows $\{a_i\}_{i=1}^t$ in $Q$.  Then
$\Lambda$ is Gorenstein if and only if $\Gamma$ is Gorenstein.
\end{cor}

\subsection{Singularity categories}

Our aim in this subsection is to show that the singularity categories
of the algebras under an arrow removal are triangle equivalent.

For this we have the following lemma in the abstract setting of cleft
extensions of abelian categories.

\begin{lem}\label{lem:singcatcleftextab}
Let $(\B,\A, \me, \ml, \mi)$ be a cleft extension of abelian
categories with enough projectives. Consider the following conditions.
\begin{enumerate}[\rm(i)]
\item $\sup\{ \pd_\B \me(P)\mid P\in\Proj(\A)\} = p_\A$ for some integer $p_\A$. 
\item $\sup\{ \pd_\A \mi(F)\mid F\in\Proj(\B)\} = p_\B$ for some integer $p_\B$. 
\item $\sup\{ \pd_\A H(B) \mid B\in\B\} = n_H$ for some integer
  $n_H$. 
\item $\sup\{ \pd_\A G(A) \mid A\in\A\} = n_G$ for some integer  $n_G$. 
\end{enumerate}
\begin{enumerate}[\rm(a)]
\item If \textup{(ii)} holds, then $\mi\colon \B\to\A$ induces a
  functor $\mi\colon \Dsg(\B) \to \Dsg(\A)$. 
\item If \textup{(i)} holds, then $\me\colon \A \to \B$ induces a
  functor $\me\colon \Dsg(\A) \to \Dsg(\B)$. 
\item If \textup{(i)} and \textup{(ii)} hold, then $\me\mi\colon
  \Dsg(\B) \to \Dsg(\B)$ is isomorphic to the identity functor. 
\item If $\ml$ is an exact functor, then $\ml\colon \B \to \A$ induces
  a functor $\ml\colon \Dsg(\B) \to \Dsg(\A)$. 
\item If \textup{(i)} and \textup{(iv)} hold and $\ml$ is an exact
  functor, then $\ml\me\colon \Dsg(\A) \to \Dsg(\A)$ is isomorphic to
  the identity functor. 
\item If \textup{(i)}--\textup{(iv)} hold and the
  functor $\ml$ is exact, then $\me\colon \Dsg(\A) \to \Dsg(\B)$ is a
  singular equivalence. 
\end{enumerate}
\end{lem}
\begin{proof}
(a) Since the functor $\mi\colon \B \to \A$ is exact, we have an
induced functor $\mi\colon \mD(\B) \to \mD(\A)$.  By property (ii) the
functor $\mi$ induce a functor $\mi\colon \mD^\perf(\B) \to
\mD^\perf(\A)$.  The claim follows from this.

(b) This follows as the claim in (a).

(c) This follows from (a) and (b) and the fact that $\me\mi \simeq
\iden_\B$. 

(d) Since the functor $\ml\colon \B \to \A$ preserves projective
objects, the claim is immediate.

(e) By (b) and (d) the functors $\me$ and $\ml$ induce functors on the
singularity categories.  Having the exact sequence
\[
0\to G(A) \to \ml\me(A) \to A \to 0
\]
from \eqref{eq:G-functor} and property (iv) ensure that the
composition of $\ml$ and $\me$ is isomorphic to the identity.

(f) By (c) the composition of $\me$ and $\mi$ is the identity functor
on $\Dsg(\B)$.  From (e) the composition of $\ml$ and $\me$ is
isomorphic to the identity functor on $\Dsg(\A)$.  Using the exact
sequence of functors 
\[
0\to H\to \ml \to \mi \to 0
\]
from \eqref{eq:H-functor} and property (iii), we infer that
$\ml\me$ and $\mi\me$ are isomorphic as endofunctors of $\Dsg(\A)$.
The claim follows from this. 
\end{proof}

As a consequence of Lemma~\ref{lem:singcatcleftextab} and
Theorem~\ref{cleftextarrowrem} we have the following. This result
constitutes part (iii) of the Main Theorem presented in the
Introduction. Below the singularity category $\Dsg(\Lambda)$ of $\Lambda$ is the Verdier quotient $\mD^b(\smod\Lambda)/\mD^\perf(\Lambda)$. 

\begin{cor}
  Let $\Lambda = kQ/I$ be an admissible quotient of a path algebra
  $kQ$ over a field $k$ and suppose that
  $\Gamma=\Lambda/\langle \{\overline{a}_i\}_{i=1}^t\rangle$ is an
  arrow removal of $\Lambda$ for the arrows $\{a_i\}_{i=1}^t$ in $Q$.
  Then the functor $\me\colon\smod\Lambda\to \smod\Gamma$ induces a
  singular equivalence between $\Lambda$ and $\Gamma\colon$
\[
\xymatrix{
\me\colon \Dsg(\Lambda) \ar[r]^{ \ \simeq} & \Dsg(\Gamma) }
\] 
\end{cor}

The next example shows that algebras can be of finite, tame or wild
representation type and still be singular equivalent to each other.

\begin{exam}
\label{examsing}
Let $Q_n$ be the quiver given by 
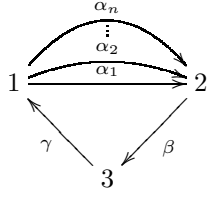
\[
\xymatrix{
1
\ar[rr]^{\alpha_1}\ar@/^0.7pc/[rr]^{\alpha_2}="a2"\ar@/^2pc/[rr]^{\alpha_n}="an"
& & 2\ar[dl]^\beta \\ 
 & 3\ar[lu]^\gamma & \ar@{..} "a2";"an"
}
\]
for $n\geq 1$.  For a field $k$ consider the relations
$\rho = \{ \alpha_1\beta, \beta\gamma, \gamma\alpha_1\}$ in
  $kQ_n$, and define the algebra
$\Lambda_n = kQ_n/\langle \rho \rangle$.  Then the algebras
  $\Lambda_1$ and $\Lambda_n$ are related by arrow removal for all
  $n \geq 2$, so that they are all singular equivalent by the above
  corollary, where $\Lambda_1$ is of finite type, $\Lambda_2$ is of
  tame type and $\Lambda_n$ is wild type for $n\geq 3$.
\end{exam}

\section{Cleft extensions and the \fg{} condition}\label{section:5}

This section is devoted to study the behaviour of the \fg{} condition
for Hochschild cohomology under the arrow removal operation. As
mentioned in the Main result of the Introduction, we prove that the
\fg{} condition is invariant under an arrow removal.  Recall
 from \cite{EHSST, Solberg:contemp} that an algebra $\Lambda$
  over a commutative ring $k$ such that $\Lambda$ is flat as a module
  over $k$ satisfies the \fg{} condition if the following is true:
\begin{enumerate}[\rm(i)]
\item The Hochschild cohomology ring $\HH^*(\Lambda)$ of $\Lambda$ is
  noetherian. 
\item The $\HH^*(\Lambda)$-module $\Ext^*_\Lambda(\Lambda/\rad
  \Lambda, \Lambda/\rad \Lambda)$ is finitely generated.
\end{enumerate}
  
Towards this we start with the following result where we show that
starting with an arrow removal and passing to the corresponding
enveloping algebras we still get a cleft extension.

\begin{prop}\label{prop:cleftenv}
  Let $\Lambda=kQ/I$ be an admissible quotient of a path algebra $kQ$,
  and suppose $\{a_i\}_{i=1}^t$ is a set set of arrows in $Q$
    such that $\Gamma = \Lambda/\langle
    \{\overline{a}_i\}_{i=1}^t\rangle$ is an arrow removal.  Let
  $\nu\colon \Gamma \to \Lambda$ and $\pi\colon \Lambda \to \Gamma$ be
  the algebra homomorphism defining the cleft extension.  Then the
  following assertions hold.
\begin{enumerate}[\rm(i)]
\item The algebra homomorphisms 
\[\nu\otimes\nu \colon \Gamma^{\opp}\otimes_k \Gamma \to
  \Lambda^{\opp}\otimes_k\Lambda\]
and 
\[\pi\otimes\pi \colon \Lambda^{\opp}\otimes_k\Lambda \to
  \Gamma^{\opp}\otimes_k \Gamma\] 
defines $\Gamma^\evl=\Gamma^{\opp}\otimes_k \Gamma$ and $\Lambda^\evl =
\Lambda^{\opp}\otimes_k\Lambda$ as a cleft extension.
\item $_{\Gamma^\evl}\Lambda^\evl$ and $\Lambda^\evl_{\Gamma^\evl}$
  are projective modules.
\item The restriction functor $\me^\evl$ along the algebra
  homomorphism $\nu\otimes\nu$ preserves projective modules (and is
  exact), and the functor 
\[\ml^\evl= - \otimes_{\Gamma^\evl}\Lambda^\evl\colon
  \smod\Gamma^\evl\to \smod\Lambda^\evl\]
is exact.  In particular the condition
\textup{(\ref{assumptionscleftext})} is satisfied for the cleft extension
$\Lambda^\evl\to \Gamma^\evl$. 
\end{enumerate}
\end{prop}
\begin{proof}
(i) It is straightforward to see that $(\pi\otimes\pi)(\nu\otimes\nu)
= \id_{\Gamma^\evl}$. 

(ii) Since $\Lambda$ and $\Gamma$ is an arrow removal,
$_\Gamma\Lambda$ and $\Lambda_\Gamma$ are projective modules.  Since
$_{\Gamma^\evl}\Lambda^\evl \simeq {_\Gamma\Lambda^{\opp}}\otimes_k
\Lambda_\Gamma$, it follows that $_{\Gamma^\evl}\Lambda^\evl$ is a
projective module over $\Gamma^\evl$.  Similarly we infer that
$\Lambda^\evl_{\Gamma^\evl}$ is a projective module over
$\Gamma^\evl$.

(iii) Both of the claims follows from (ii).
\end{proof}

The functors $F$ and $G$ are crucial for a cleft extension.  Next we
see how the $F$- and the $G$-functors are connected for a cleft
extension of algebras and the corresponding cleft extension for the
enveloping algebras.

\begin{lem}\label{lem:FG-properties}
  Let $\Lambda$ and $\Gamma$ be a cleft extension given by the algebra
  homomorphisms
  $\Gamma\xrightarrow{\nu} \Lambda \xrightarrow{\pi} \Gamma$. Then for
  a $\Gamma$-bimodule $B$ the following hold.
\begin{enumerate}[\rm(a)]
\item The endofunctor $F$ of $\smod\Gamma$
  applied to $B$ defines a $\Gamma$-bimodule and the exact sequence
\[
0\to F(B) \to B\otimes_\Gamma \Lambda
  \xrightarrow{\operatorname{mult}(1\otimes\pi)} B\to 0
\]
obtained from \eqref{splitsequencegivesF} splits as a sequence
of $\Gamma$-bimodules.  
\item We have the isomorphism 
\[
F(\Lambda\otimes_\Gamma B) \simeq \Lambda\otimes_\Gamma F(B).
\]
\item Let $F^{\opp}$ be the endofunctor of $\smod\Gamma^{\opp}$ considering
  $\Gamma^{\opp}$ and $\Lambda^{\opp}$ as a cleft extension of
  algebras.  We have  
\[
F^\evl(B)  \simeq (\Lambda\otimes_\Gamma F(B_\Gamma)) \oplus F^{\opp}({_\Gamma B}).
\]
\end{enumerate}
For a $\Lambda$-bimodule $B$ the following hold. 
\begin{enumerate}[\rm(a)]
\setcounter{enumi}{3}
\item The endofunctor $G$ of
  $\smod\Lambda$ applied to $B$ defines a $\Lambda$-bimodule. 
\item When $\Lambda_\Gamma$ is projective, we have
\[
G(\Lambda\otimes_\Gamma (_\Gamma B_\Lambda)) = \Lambda\otimes_\Gamma G(B).
\]
\item Let $G^{\opp}$ be the endofunctor of $\smod\Lambda^{\opp}$ considering
$\Gamma^{\opp}$ and $\Lambda^{\opp}$ as a cleft extension of algebras.
  For a $\Lambda$-bimodule $B$ we have an exact sequence
\[
0\to \Lambda\otimes_\Gamma G(B_\Lambda) \to G^\evl(B) \to
  G^{\opp}({_\Lambda B})\to 0, 
\]
when $\Lambda_\Gamma$ is projective.
\item In this final statement let $\Lambda$ and $\Gamma$ be an arrow removal given by a set of
  arrows $\{\overline{a}_i\}_{i=1}^t$ in the quiver of the algebra
  $\Lambda$ as defined in subsection
  \textup{\ref{subsec:cleftextfromarrowremoval}. 
  Then the following hold.}
\begin{enumerate}[\rm(i)]
\item $(F^\evl)^2(\me^\evl(\Lambda)) = 0$, 
\item $(G^\evl)^2(\Lambda) = 0$, 
\item $G^\evl(\Lambda)$ is a projective $\Lambda$-bimodule,
\end{enumerate}
\end{enumerate}
\end{lem}
\begin{proof}
(a) Let $B$ be a $\Gamma$-bimodule and consider the exact sequence
\[0\to F(B_\Gamma) \to B\otimes_\Gamma \Lambda
  \xrightarrow{\operatorname{mult}(1\otimes\pi)} B_\Gamma \to 0,\]
where the map $\textrm{mult}$ is a homomorphism of
$\Gamma$-bimodules. This implies that $F(B_\Gamma)$ is a
$\Gamma$-bimodule whenever $B$ is a $\Gamma$-bimodule.  The above
exact sequence splits as right $\Gamma$-modules by Lemma
\ref{lemsplitexactseq}, but the splitting $(1\otimes\nu)
\textrm{mult}^{-1}$ is also a homomorphism of
$\Gamma$-bimodules.  Hence the final claim follows. 

(b) Let $B$ be a $\Gamma$-bimodule and consider the exact
sequence
\[0\to F(B_\Gamma) \to B\otimes_\Gamma \Lambda_\Gamma 
  \xrightarrow{\textrm{mult}(1\otimes \pi)}  B_\Gamma \to 0,\]
which splits as an exact sequence of $\Gamma$-bimodules.  Tensoring
this split exact sequence with $\Lambda\otimes_\Gamma-$ we get the 
following exact commutative diagram
\[\xymatrix{ 
0 \ar[r] & \Lambda\otimes_\Gamma F(B_\Gamma)\ar[r]\ar[d]^\simeq
            &  \Lambda\otimes_\Gamma B\otimes_\Gamma \Lambda_\Gamma \ar[rr]^-{1\otimes
    \textrm{mult}(1\otimes \pi)} \ar@{=}[d] & & \Lambda\otimes_\Gamma B_\Gamma
  \ar[r]\ar@{=}[d]  & 0\\
0 \ar[r] & F(\Lambda\otimes_\Gamma B_\Gamma)\ar[r] 
            &  \Lambda\otimes_\Gamma B\otimes_\Gamma \Lambda_\Gamma
            \ar[rr]^-{\textrm{mult}((1_{\Lambda\otimes_\Gamma B})\otimes \pi)} & & \Lambda\otimes_\Gamma B_\Gamma
  \ar[r] & 0
}\]
The claim follows from this. 

(c) Recall that $F^\evl$ is given by the exact sequence 
\[0\to F^\evl \to \me^\evl\ml^\evl \to \iden_{\smod\Gamma^\evl} \to 0\]
Let $B$ be a $\Gamma$-bimodule. Then $\ml^\evl(B) =
\Lambda\otimes_\Gamma B\otimes_\Gamma \Lambda$, so that
\[\me^\evl\ml^\evl(B) = {_\Gamma\Lambda\otimes_\Gamma B\otimes_\Gamma\Lambda_\Gamma}.\] 
We construct the following commutative diagram 
\[\xymatrix{
 &  & & & & 0\ar[d] & \\
 & 0\ar[d] & 0\ar[d] & & & F^{\opp}({_\Gamma B})\ar[d] & \\
0\ar[r] & \Lambda\otimes_\Gamma F(B_\Gamma)\ar[r]\ar[d] & \Lambda\otimes_\Gamma
B\otimes_\Gamma \Lambda \ar[rrr]^-{1\otimes \textrm{mult}(1\otimes\pi)}\ar@{=}[d]
& & & \Lambda\otimes_\Gamma B \ar[r]\ar[d]^{\textrm{mult}(\pi\otimes 1)} & 0\\
0\ar[r] & F^\evl(B)\ar[d]\ar[r] & \Lambda\otimes_\Gamma
B\otimes_\Gamma \Lambda \ar[rrr]^-{(\textrm{mult}(\pi\otimes 1))(1\otimes (\textrm{mult}(1\otimes\pi)))}\ar[d]
& & & B \ar[r]\ar[d] & 0\\
& F^{\opp}({_\Gamma B})\ar[d] & 0 & & & 0 &\\
& 0 &  & & & &\\
}\]
where the second row is split exact by Lemma \ref{lemsplitexactseq}.
This implies the first isomorphism below
\begin{align}
\Lambda\otimes_\Gamma B\otimes_\Gamma \Lambda & \simeq F^\evl(B)\oplus
                                                B\notag\\
  &  \simeq \Lambda\otimes_\Gamma F(B)\oplus F^{\opp}(B) \oplus B\notag
\end{align}
Since the first row in the above diagram is a split exact sequence by
(a) and the first column is a pullback of the first row, the second
isomorphism follows.  Cancelling the direct summand $B$ on each side
implies that
$F^\evl(B) \simeq \Lambda\otimes_\Gamma F(B)\oplus F^{\opp}(B)$.

(d) Let $B$ be a $\Lambda$-bimodule and consider the exact sequence
\begin{equation}\label{eq:G-sequence}
0\to G(B_\Lambda) \to B\otimes_\Gamma \Lambda
  \xrightarrow{\textrm{mult}} B_\Lambda \to 0,
\end{equation}
where the map $\textrm{mult}$ is a homomorphism of
$\Lambda$-bimodules. This implies that $G(B_\Lambda)$ is a
$\Lambda$-bimodule whenever $B$ is a $\Lambda$-bimodule. 

(e) Let $B$ be a $\Lambda$-bimodule. Since $\Lambda_\Gamma$ is
projective, tensoring the exact sequence \eqref{eq:G-sequence} with
$\Lambda\otimes_\Gamma -$ leaves it exact and we obtain the following
commutative diagram 
\[\xymatrix{ 
0 \ar[r] & \Lambda\otimes_\Gamma G(B_\Lambda)\ar[r]\ar@{=}[d] 
            &  \Lambda\otimes_\Gamma B\otimes_\Gamma \Lambda
            \ar[rr]^-{1_\Lambda\otimes\textrm{mult}_B} \ar@{=}[d] & &
            \Lambda\otimes_\Gamma B \ar[r]\ar@{=}[d]  & 0\\
0 \ar[r] & G(\Lambda\otimes_\Gamma B)\ar[r] 
            &  \Lambda\otimes_\Gamma B\otimes_\Gamma \Lambda
            \ar[rr]^-{\textrm{mult}_{\Lambda\otimes_\Gamma B}} & & \Lambda\otimes_\Gamma B
  \ar[r] & 0
}\]
The claim follows from this. 

(f) This follows in a similar way as for $F^\evl$, and it left to the
reader. 

(g) Let  $\Lambda$ and $\Gamma$ be an arrow removal, and let $B$ be a
$\Gamma$-bimodule. Then 
\begin{align}
(F^\evl)^2(B) & = F^\evl(F^\evl(B)),\notag\\
& \simeq F^\evl(\Lambda\otimes_\Gamma F(B) \oplus F^{\opp}(B)), \textrm{\quad using (c)}\notag\\
& \simeq F^\evl(\Lambda\otimes_\Gamma F(B)) \oplus
  F^\evl(F^{\opp}(B)), \textrm{\ using additivity}\notag\\
& = \Lambda\otimes_\Gamma F((\Lambda\otimes_\Gamma F(B_\Gamma))_\Gamma) \oplus
  F^{\opp}({_\Gamma(\Lambda\otimes_\Gamma F(B_\Gamma))})\notag\\
& \qquad\qquad \oplus \Lambda\otimes_\Gamma F(F^{\opp}({_\Gamma B})_\Gamma)
  \oplus F^{\opp}(F^{\opp}({_\Gamma B})), \textrm{\quad using (c)}\notag\\
& = \Lambda\otimes_\Gamma F^2(\Lambda\otimes_\Gamma B_\Gamma) \oplus
  F^{\opp}F(\Lambda\otimes_\Gamma B)\notag\\
 & \qquad\qquad \oplus \Lambda\otimes_\Gamma F(F^{\opp}({_\Gamma
   B})_\Gamma) \oplus F^{\opp}(F^{\opp}(B)), \textrm{\quad using (b)}\notag
\end{align}
Since $F^2 = 0$ and $(F^{\opp})^2 = 0$ for an arrow removal by Theorem
\ref{propremoveiscleft} (ii) (g), we have 
\[(F^\evl)^2(B) = F^{\opp}F(\Lambda\otimes_\Gamma B) \oplus
  \Lambda\otimes_\Gamma FF^{\opp}(B).\]
When we let $\langle\{\overline{a}_i\}_{i=1}^t\rangle$ denote the
$\Gamma$-sub-bimodule of $\Lambda$ generated by
$\{\overline{a}_i\}_{i=1}^t$,  we have by Lemma \ref{lem:F-func-prop}
that
\[F^{\opp}F(\Lambda\otimes_\Gamma B) = \langle\{\overline{a}_i\}_{i=1}^t\rangle
  \otimes_\Gamma \Lambda\otimes_\Gamma B\otimes_\Gamma
  \langle\{\overline{a}_i\}_{i=1}^t\rangle\]
 and  
\[FF^{\opp}(B) = \langle\{\overline{a}_i\}_{i=1}^t\rangle \otimes_\Gamma B
  \otimes_\Gamma \langle\{\overline{a}_i\}_{i=1}^t\rangle.\]
When we specialize to
$B = {_\Gamma\Lambda_\Gamma} = \me^\evl(\Lambda)$ and use that
$\Lambda \simeq \Gamma
\oplus\langle\{\overline{a}_i\}_{i=1}^t\rangle$, then
\begin{align}
FF^{\opp}(\me^\evl(\Lambda)) & =  \langle\{ \overline{a}_i \}_{i=1}^t\rangle
  \otimes_\Gamma (\Gamma\oplus \langle\{ \overline{a}_i
  \}_{i=1}^t\rangle)\otimes_\Gamma \langle\{ \overline{a}_i \}_{i=1}^t\rangle\notag\\
& \simeq \langle\{ \overline{a}_i \}_{i=1}^t\rangle\otimes_\Gamma \langle\{ \overline{a}_i
  \}_{i=1}^t\rangle\notag\\
& \qquad\qquad \oplus \langle\{ \overline{a}_i \}_{i=1}^t\rangle\otimes_\Gamma \langle\{
  \overline{a}_i \}_{i=1}^t\rangle\otimes_\Gamma \langle\{
  \overline{a}_i \}_{i=1}^t\rangle\notag\\
& = 0,\notag
\end{align}
since $f_j\Gamma e_i = 0$ for all $i,j = 1,2,\ldots, t$.  For similar
reasons we obtain that $F^{\opp}F(\Lambda\otimes_\Gamma \Lambda) = 0$
and consequently 
\[(F^\evl)^2(\me^\evl(\Lambda)) = 0.\]
Since $\me^\evl(G^\evl(B))^2 \simeq (F^\evl)^2(\me^\evl(B))$ by
  Lemma \ref{lemnilpotentfunctors} (i) and $\me^\evl$ is faithful, we
infer that $(G^\evl)^2(\Lambda) = 0$.  Using similar arguments as
above
$F^\evl(\me^\evl(\Lambda)) \simeq \langle\{ \overline{a}_i
\}_{i=1}^t\rangle^{\oplus 2}$ as a $\Gamma$-bimodule.  For an arrow
removal $\langle\{ \overline{a}_i \}_{i=1}^t\rangle$ is a projective
$\Gamma$-bimodule.  Then by Lemma \ref{lemnilpotentfunctors} (i)
$\me^\evl(G^\evl(\Lambda)) \simeq F^\evl(\me^\evl(\Lambda))$ and it is
projective.  Since the functor $\ml^\evl$ preserves projective
modules, the bimodule $\ml^\evl\me^\evl G^\evl(\Lambda)$ is
projective.  We have the exact sequence
\[0\to (G^\evl)^2(\Lambda) \to \ml^\evl\me^\evl G^\evl(\Lambda) \to
  G^\evl(\Lambda) \to 0,\]
which implies that $G^\evl(\Lambda) \simeq \ml^\evl\me^\evl
G^\evl(\Lambda)$ is a projective $\Lambda$-bimodule.  
\end{proof}

The following result establishes a close relationship between the
Hochschild cohomology rings for the algebras in an arrow removal. {The
  interested reader is suggested to compare the isomorphism below with
  \cite[Theorem~4.6]{CLMS}.

\begin{prop}
 If $\pi\colon\Lambda\to\Gamma$ is an arrow removal, then
\[\Ext^*_{\Lambda^\evl}(\Lambda,\Lambda) \simeq
\Ext^*_{\Gamma^\evl}(\Gamma,\Gamma \oplus \Ker\pi)\]
is an isomorphism for $* > 1$.
\end{prop}
\begin{proof}
As above we have the exact sequence 
\[0\to G^\evl(\Lambda) \to \ml^\evl\me^\evl(\Lambda) \to \Lambda \to
  0.\] 
By Lemma \ref{lem:FG-properties} (g) the bimodule $G^\evl(\Lambda)$ is
projective.  The condition (\ref{assumptionscleftext}) is satisfied
for the cleft extension $\Lambda^\evl \to \Gamma^\evl$ (see
Proposition \ref{prop:cleftenv} (iii)), so that we can
use Lemma~\ref{lem:commdiagcleft} to obtain 
\[\Ext^*_{\Lambda^\evl}(\Lambda,\Lambda) \simeq
  \Ext^*_{\Gamma^\evl}(\me^\evl(\Lambda),\me^\evl(\Lambda))\] 
for $* > 1$.  The restriction
$\me^\evl(\Lambda) \simeq \Gamma \oplus \Ker\pi$, 
where $\Ker\pi = \langle \{\overline{a}_i\}_{i=1}^t\rangle$  is a
projective $\Gamma$-bimodule.  This implies that
\[\Ext^*_{\Lambda^\evl}(\Lambda,\Lambda) \simeq
\Ext^*_{\Gamma^\evl}(\Gamma,\Gamma \oplus \Ker\pi),\]
for $* > 1$ and it completes the proof. 
\end{proof}

For the \fg{}-property to be preserved for an arrow removal, not only
the Hochschild cohomology rings need to be related, but also their
action on the $\Ext$-groups must respect each other, in order to apply
general results from \cite[Proposition 6.4]{PSS}.  The following two
results prepares for this.
\begin{lem}
Let $\pi\colon\Lambda\to\Gamma$ be an arrow removal.
Let $M$ be a right $\Lambda$-module and $B$ a $\Lambda$-bimodule. Then
the map 
\[\me(M)\otimes_\Gamma\me^\evl(B) \xrightarrow{\varphi}
  \me(M\otimes_\Lambda B)\]
given by $m\otimes b \mapsto m\otimes b$ is well-defined, functorial
in both variables, and an onto map of right $\Gamma$-modules. 
\end{lem}
\begin{proof}
The module $\me(M)\otimes_\Gamma\me^\evl(B) = M_\Gamma \otimes_\Gamma
{_\Gamma B_\Gamma}$ and the module $\me(M\otimes_\Lambda B) =
M\otimes_\Lambda B_\Gamma$. Therefore the map $\varphi$ is the natural
projection. 
\end{proof}

\begin{prop}\label{prop:commutativeFGdiagram}
Let $\pi\colon\Lambda\to\Gamma$ be an arrow removal.
The following diagram is commutative
\[\xymatrix{
\Ext^*_{\Lambda^\evl}(\Lambda,\Lambda) \ar[rr]^{M\otimes_\Lambda-}
\ar[ddd]_{\me^\evl} & & \Ext^*_\Lambda(M,M)\ar[d]^\me \\
 & & \Ext^*_\Gamma( \me(M),\me(M) )\ar[d]^{\Ext^*_\Gamma(\varphi,-)} \\
 & & \Ext^*_\Gamma( \me(M)\otimes_\Gamma\me^\evl(\Lambda), \me(M))\\
\Ext^*_{\Gamma^\evl}(\me^\evl(\Lambda), \me^\evl(\Lambda))
\ar[rr]^-{\me(M)\otimes_\Gamma-} & & \Ext^*_\Gamma(
\me(M)\otimes_\Gamma\me^\evl(\Lambda), \me(M)\otimes_\Gamma
\me^\evl(\Lambda)) \ar[u]_{\Ext^*_\Gamma(-,\varphi)}
}\]
\end{prop}
\begin{proof}
Let $\eta\colon \Omega^n_{\Lambda^\evl}(\Lambda)\to \Lambda$ represent
an element in $\Ext^n_{\Lambda^\evl}(\Lambda,\Lambda)$. As an
extension $\eta$ correspond to the lower row in the following
commutative diagram
\[\xymatrix{
0\ar[r] & \Omega^n_{\Lambda^\evl}(\Lambda) \ar[r]\ar[d]_\eta &
P_{n-1}\ar[r]\ar[d] & P_{n-2}\ar[r]\ar@{=}[d] & \cdots\ar[r] &
P_0\ar[r]\ar@{=}[d] & \Lambda\ar[r] \ar@{=}[d] & 0\\
0\ar[r] & \Lambda\ar[r] & E\ar[r] & P_{n-2}\ar[r] & \cdots\ar[r] &
P_0\ar[r] & \Lambda\ar[r] & 0
}\]
where the first row is the start of a projective resolution of
$\Lambda$ over $\Lambda^\evl$.  Tensoring this diagram with $M$ over
$\Lambda$ we obtain the extension $M\otimes_\Lambda\eta$ as the lower
row in the following exact commutative diagram
\[
\xymatrix@C7pt{
0\ar[r] & M\otimes_\Lambda\Omega^n_{\Lambda^\evl}(\Lambda) \ar[r]\ar[d]_{M\otimes\eta} &
M\otimes_\Lambda P_{n-1}\ar[r]\ar[d] & M\otimes_\Lambda P_{n-2}\ar[r]\ar@{=}[d] & \cdots\ar[r] &
M\otimes_\Lambda P_0\ar[r]\ar@{=}[d] & M\otimes_\Lambda\Lambda\ar[r] \ar@{=}[d] & 0\\
0\ar[r] & M\otimes_\Lambda\Lambda\ar[r] & M\otimes_\Lambda E\ar[r] & M\otimes_\Lambda P_{n-2}\ar[r] & \cdots\ar[r] &
M\otimes_\Lambda P_0\ar[r] & M\otimes_\Lambda\Lambda\ar[r] & 0
}
\]
Restricting all the homomorphisms and all the modules to $\Gamma$ in
the above diagram we obtain the extension
$\me(M\otimes_\Lambda \eta)$. We use similar arguments as in the proof
of Lemma~\ref{lem:commdiagcleft} to construct it.  We first look
  at the case $n=1$ to illustrate this.  In the following commutative
  diagram, the second row is the image in
  $\Ext^1_\Gamma(\me(M),\me(M))$ and the third row is the image in
  $\Ext^1_\Gamma(\me(M)\otimes_\Gamma\me^\evl(\Lambda),
  \me(M)\otimes_\Gamma \me^\evl(\Lambda))$.
\[
\xymatrix@C7pt{
0\ar[r] & \me(M\otimes_\Lambda\Omega^1_{\Lambda^\evl}(\Lambda))
\ar[r]\ar[d]_{\me(M\otimes\eta)} & \me(M\otimes_\Lambda P_0) \ar[r]\ar[d] &
\me(M\otimes_\Lambda \Lambda) \ar[r]\ar@{=}[d] & 0\\
0\ar[r] & \me(M\otimes_\Lambda\Lambda)\ar[r] & \me(M\otimes_\Lambda E)
\ar[r] & \me(M\otimes_\Lambda \Lambda) \ar[r] & 0\\
0\ar[r] & \me(M)\otimes_\Gamma \me^\evl(\Lambda) \ar[r]\ar[u]^\varphi
& \me(M)\otimes_\Gamma \me^\evl(E) \ar[r]\ar[u]^\varphi
& \me(M)\otimes_\Gamma \me^\evl(\Lambda) \ar[r]\ar[u]^\varphi & 0
} \]
Then the pullback of the second row along $\varphi$ is equivalent to
the pushout of the third row along $\varphi$, which shows the claim
for $n = 1$.  For $n > 1$ we have the following. 
\[
\xymatrix@C7pt{
0\ar[r] & M\otimes_\Lambda\Omega^n_{\Lambda^\evl}(\Lambda)_\Gamma
\ar[r]\ar[d]_{M\otimes\eta} & 
M\otimes_\Lambda P_{n-1}\ar[r]\ar[d] & M\otimes_\Lambda P_{n-2}\ar[r]\ar@{=}[d] & \cdots\ar[r] &
M\otimes_\Lambda P_1\ar[r]\ar@{=}[d] & M\otimes_\Lambda P_0\ar[r]\ar@{=}[d] &
M\otimes_\Lambda\Lambda_\Gamma\ar[r] \ar@{=}[d] & 0\\ 
0\ar[r] & M\otimes_\Lambda\Lambda_\Gamma\ar[r] & M\otimes_\Lambda
E\ar[r] & M\otimes_\Lambda P_{n-2}\ar[r] & \cdots\ar[r] &
M\otimes_\Lambda P_1\ar[r] &
M\otimes_\Lambda P_0\ar[r] & M\otimes_\Lambda\Lambda_\Gamma\ar[r] &
0\\
0\ar[r] & M\otimes_\Lambda\Lambda_\Gamma\ar[r]\ar@{=}[u] & M\otimes_\Lambda
E\ar[r]\ar@{=}[u] & M\otimes_\Lambda P_{n-2}\ar[r]\ar@{=}[u] &
\cdots\ar[r] & M\otimes_\Lambda P_1\ar[r]\ar@{=}[u] &
E''\ar[r]\ar[u] & M\otimes_\Gamma\Lambda_\Gamma\ar[r]\ar[u]^\varphi &
0\\
0\ar[r] & M\otimes_\Lambda\Lambda_\Gamma\ar[r]\ar@{=}[u] & 
E'\ar[r]\ar[u] & M\otimes_\Gamma P_{n-2}\ar[r]\ar[u]^\varphi & \cdots\ar[r] &
M\otimes_\Gamma P_1\ar[r]\ar[u]^\varphi &
M\otimes_\Gamma P_0\ar[r]\ar[u] & M\otimes_\Gamma\Lambda_\Gamma\ar[r]\ar@{=}[u] &
0\\
0\ar[r] & M\otimes_\Gamma\Lambda_\Gamma\ar[r]\ar[u]^\varphi & M\otimes_\Gamma
E\ar[r]\ar[u] & M\otimes_\Gamma P_{n-2}\ar[r]\ar@{=}[u] & \cdots\ar[r]
& M\otimes_\Lambda P_1\ar[r]\ar@{=}[u] &
M\otimes_\Gamma P_0\ar[r]\ar@{=}[u] & M\otimes_\Gamma\Lambda_\Gamma\ar[r]\ar@{=}[u] & 0
}
\]
As said above, the second row is $\me(M\otimes_\Lambda\eta)$, the
third row is 
\[\Ext^n_\Gamma(\varphi,-)(\me(M\otimes_\Lambda \eta)),\] 
the fifth row is $\me(M)\otimes_\Gamma \me^\evl(\eta)$, and
the fourth row is 
\[\Ext^n_\Gamma(-,\varphi)(\me(M)\otimes_\Gamma \me^\evl(\eta)).\]  
The diagram shows that the extension on the third row and the
extension on the fourth row are equivalent.  In other words, the
diagram in the proposition is commutative.
\end{proof}

Next we prove the main result of this section which shows that the
\fg{} condition is invariant under the arrow removal operation. This
is part (iv) of the Main Theorem presented in the Introduction.

\begin{thm}
\label{thmFg}
Let $\Lambda=kQ/I$ be an admissible quotient of a path algebra $kQ$,
and suppose that $\Lambda\to \Gamma$ is an arrow removal.  Then
  $\Lambda$ satisfies \fg{} if and only if $\Gamma$ satisfies \fg{}.
\end{thm}
\begin{proof} 
  We use \cite[Proposition 6.4]{PSS} with $N = N' = \Gamma/\rad\Gamma$
  and $M = M' = \mi(\Gamma/\rad\Gamma)$, where $N$ is the direct sum
  of all simple $\Gamma$-modules and $M$ is the direct sum of all
  simple $\Lambda$-modules.

We have that
\begin{align}
\me(M)\otimes_\Gamma \me^\evl(\Lambda) & \simeq
\me(M)\otimes_\Gamma (\Gamma \oplus \Ker\pi)\notag\\
& \simeq \me(M)\oplus (\me(M)\otimes_\Gamma \Ker\pi)\notag\\
& \simeq \me(M) \oplus F(\me(M))\notag
\end{align}
Since $F(\me(M))$ is projective by Theorem
\ref{cleftextarrowrem} (f), the homomorphism 
$\Ext^*_\Gamma(\varphi,-)$ is an isomorphism for $* > 0$ in the
commutative diagram in Proposition
\ref{prop:commutativeFGdiagram}. Since $\Lambda$ is Gorenstein if and
only if $\Gamma$ is Gorenstein by Theorem~\ref{thmGorenstein}, we
have that both $\Lambda$ and $\Gamma$ are Gorenstein whenever we
assume one of them is Gorenstein.  Hence if we assume that one of
$\Lambda$ and $\Gamma$ has \fg{}, then $\Gamma$ is Gorenstein by
\cite[Prop.\ 2.2]{EHSST}.  So we
can suppose $\Gamma$ is Gorenstein.  Then
$\Ker\pi$ has finite injective dimension as a
right $\Gamma$-module, say $n$.  This implies that the homomorphism
$\Ext^*_\Gamma(-,\varphi)$ is an isomorphism for $* > n$ in the
commutative of Proposition
\ref{prop:commutativeFGdiagram}.  Furthermore, $\Gamma^\evl$ and
$\Lambda^\evl$ are both also Gorenstein.
Suppose that $\Gamma^\evl$ has Gorenstein dimension $d$.

Let $p_\Gamma = \pi \colon \me^\evl(\Lambda) \to \Gamma$ be the natural
projection.  Then construct the following commutative diagram for $* >
\max\{ n, d\}$.  The upper square is the the commutative square of Proposition
\ref{prop:commutativeFGdiagram}. 
\[\xymatrix{
\Ext^*_{\Lambda^\evl}(\Lambda,\Lambda) \ar[rr]^{M\otimes_\Lambda-}
\ar[d]_{\me^\evl} & & \Ext^*_\Lambda(M,M)\ar[d]^\simeq \\
\Ext^*_{\Gamma^\evl}(\me^\evl(\Lambda), \me^\evl(\Lambda))
\ar[rr]^-{\me(M)\otimes_\Gamma-}\ar[d]_{\Ext^*_{\Gamma^\evl}(-, p_\Gamma)} & & \Ext^*_\Gamma(
\me(M)\otimes_\Gamma\me^\evl(\Lambda), \me(M)\otimes_\Gamma
\me^\evl(\Lambda)) \ar[d]^{\Ext^*_\Gamma(-,1\otimes p_\Gamma)}\\
\Ext^*_{\Gamma^\evl}(\me^\evl(\Lambda), \Gamma)
\ar[rr]^-{\me(M)\otimes_\Gamma-} & & 
\Ext^*_\Gamma(\me(M)\otimes_\Gamma\me^\evl(\Lambda),
\me(M)\otimes_\Gamma \Gamma) \\
\Ext^*_{\Gamma^\evl}(\Gamma, \Gamma)
\ar[rr]^-{\me(M)\otimes_\Gamma-}\ar[u]^{\Ext^*_\Gamma(p_\Gamma,-)}\ar@{=}[d] & & 
\Ext^*_\Gamma(\me(M)\otimes_\Gamma\Gamma,
\me(M)\otimes_\Gamma \Gamma) \ar[u]_{\Ext^*_\Gamma(1\otimes
  p_\Gamma,-)}\ar[d]^\simeq\\
\Ext^*_{\Gamma^\evl}(\Gamma, \Gamma)
\ar[rr]^-{\me(M)\otimes_\Gamma-} & &  \Ext^*_\Gamma(\me(M), \me(M))
}\]
All the vertical maps in this diagram are isomorphisms and the diagram
is commutative.  Then using \cite[Proposition 6.4]{PSS} with
$N = N' = \Gamma/\rad\Gamma$ and $M = M' = \mi(N)$,
where $M$ is
the direct sum of all simple $\Lambda$-modules and noting that
  $\me(M) = \me\mi(N) \simeq N$, we obtain that
$\Lambda$ has \fg{} if and only if $\Gamma$ has \fg{}.
\end{proof}

\begin{exam}
Let $\Lambda_n = kQ_n/\langle \rho \rangle$ be the algebra of Example~\ref{examsing}. After removing the arrows $\alpha_2, \ldots, \alpha_n$, we obtain a radical square zero Nakayama algebra which satisfies \fg{} by \cite[Proposition~1.4]{ES}. By Theorem~\ref{thmFg}, we infer that $\Lambda_n$ satisfies \fg{}.
\end{exam}

We end the paper with an example showing that a general arrow removal
(factoring out an arrow) and preserving \fg{} is not possible.

\begin{exam}
Consider the following example presented by Fei Xu \cite[3.1 The
category $\mathcal{E}_0$]{FeiXu}.  Let $Q$ be the quiver given by
\[\xymatrix{
1\ar@(ul,ur)^a\ar@(dl,dr)_b \ar[r]^c & 2
}\]
and the ideal $I = \langle a^2, ab - ba, b^2, ac\rangle$ in $kQ$ for a
field $k$.  Denote by $\Lambda$ the factor algebra $kQ/I$. By a
result in a forthcoming paper or by direct computations,
$\Lambda\xrightarrow{\pi} \Gamma = \Lambda/\langle c \rangle$ is a
cleft extension. Then by \cite{FeiXu} $\Lambda$ does not satisfy
\fg{}, while $\Gamma$ do satisfy \fg{} (since $\Gamma$ is a
symmetric radical cube zero algebra satisfying \fg{} by \cite{ES}).
\end{exam}

\end{document}